\documentclass{article}

\usepackage{epsf,latexsym,amssymb}
\usepackage{graphicx}
\usepackage{times}
\usepackage{amsmath}  

\newcommand{\comment}[1]{}
\newcommand{\edo}{\end{document}}

\newcommand{\R}{{\mathbb R}}  

\newtheorem{theorem}{Theorem}
\newtheorem{itlemma}{Lemma}[section] 
\newtheorem{itproposition}[itlemma]{Proposition}
\newtheorem{itcorollary}[itlemma]{Corollary}
\newtheorem{itremark}[itlemma]{Remark}
\newtheorem{itdefinition}[itlemma]{Definition}
\newtheorem{itexample}[itlemma]{Example}

\newenvironment{lemma}{\begin{itlemma}\rm}{\end{itlemma}} 
\newenvironment{remark}{\begin{itremark}\rm}{\end{itremark}} 
\newenvironment{corollary}{\begin{itcorollary}\rm}{\end{itcorollary}}
\newenvironment{proposition}{\begin{itproposition}\rm}{\end{itproposition}}
\newenvironment{definition}{\begin{itdefinition}\rm}{\end{itdefinition}}
\newenvironment{example}{\begin{itexample}\rm}{\end{itexample}}

\newcommand{\be}[1]{\begin{equation}\label{#1}}
\newcommand{\ee}{\end{equation}}
\newcommand{\bl}[1]{\begin{lemma}\label{#1}}
\newcommand{\ble}[1]{\begin{lemmaex}\label{#1}}
\newcommand{\br}[1]{\begin{remark}\label{#1}}
\newcommand{\bt}[1]{\begin{theorem}\label{#1}}
\newcommand{\bd}[1]{\begin{definition}\label{#1}}
\newcommand{\bp}[1]{\begin{proposition}\label{#1}}
\newcommand{\bc}[1]{\begin{corollary}\label{#1}}
\newcommand{\bex}[1]{\begin{example}\label{#1}}
\newcommand{\ec}{\mybox\end{corollary}}
\newcommand{\eex}{\mybox\end{example}}
\newcommand{\eem}{\mybox\end{example}}
\newcommand{\el}{\mybox\end{lemma}}
\newcommand{\er}{\mybox\end{remark}}
\newcommand{\et}{\qed\end{theorem}}
\newcommand{\ed}{\mybox\end{definition}}
\newcommand{\ep}{\mybox\end{proposition}}
\newcommand{\epr}{\end{proof}}
\newcommand{\bpr}{\begin{proof}}

\newcommand{\ecs}{\end{corollary}}
\newcommand{\eexs}{\end{example}}
\newcommand{\els}{\end{lemma}}
\newcommand{\ers}{\end{remark}}
\newcommand{\ets}{\end{theorem}}
\newcommand{\eds}{\end{definition}}
\newcommand{\eps}{\end{proposition}}
\newcommand{\halmos}{\rule{1ex}{1.4ex}}
\newcommand{\qed}{\hfill \halmos} 
\newcommand{\mybox}{\hfill $\Box$} 
\newcommand{\beq}{\begin{eqnarray}}
\newcommand{\eeq}{\end{eqnarray}}
\newcommand{\beqn}{\begin{eqnarray*}}
\newcommand{\eeqn}{\end{eqnarray*}}
\newcommand{\bi}{\begin{itemize}}
\newcommand{\ei}{\end{itemize}}
\newcommand{\ben}{\begin{enumerate}}
\newcommand{\een}{\end{enumerate}}

\newcommand{\bes}[1]{\begin{subequations}\label{#1}\begin{eqnarray}}
\newcommand{\ees}[1]{\end{eqnarray}\end{subequations}}

\newcommand{\st}{\, | \,}

\newenvironment{proof}{\noindent {\em Proof}.\ }{\hspace*{\fill}$\halmos$\medskip}

\newcommand{\mypmatrix}[1]{\left(\begin{array}{cccccccccccc}#1\end{array}\right)}

\title{A remark about polynomials with specified local minima\\ and no other critical points}
\author{Eduardo D. Sontag}

\begin{document}
\maketitle

\section{Introduction}

This fact must surely be well-known, but it seems worth giving a simple and
quite explicit proof:

\bp{prop:main}
Take any finite subset $X$ of $\R^n$, $n\geq 2$.
Then, there is a polynomial function $P:\R^n\rightarrow \R$ which has local minima
on the set $X$, and has no other critical points.
\eps

A weaker version, constructing one particular subset $X$ of $\R^2$ with the
stated property, was shown in a beautiful little note by an undergraduate, Ian
Robertson, as part of an REU conducted by Alan Durfee \cite{robertson92};
references to this work, and the context from degree theory, can be found
in~\cite{durfee98}.  Our construction generalizes Robertson's.

We were interested in this question because of the following consequence,
obtained immediately by considering the gradient vector field
$g(x)=-\nabla P(x)$:

\bc{cor:main}
Take any finite subset $X$ of $\R^n$, $n\geq 2$.
Then, there is a polynomial vector field $g:\R^n\rightarrow \R^n$ which has
asymptotically stable equilibria on the set $X$, and has no other equilibria.
\ecs

Suppose that $g$ is as in Corollary~\ref{cor:main} and that $X$ has more than
one element.  Let ${\cal O}_x$, $x$$\in $$X$, be the domains of attraction of the points
in $X$.  The union of the sets ${\cal O}_x$ cannot equal all of $\R^n$, since these
are disjoint open sets; pick any $\xi \in \R^n \setminus \bigcup _x {\cal O}_x$.
The omega-limit set $\Omega ^+(\xi )$ cannot intersect $X$ (since points of $X$ are
asymptotically stable).  Thus there are points that do not converge to any
equilibria.  (Alternatively, one could arrive at the same conclusion appealing
to topological degree arguments.)  Thus, the following is also of interest
(and much easier to prove).

\bp{prop:weaker}
Take any finite subset $X$ of $\R^n$, $n\geq 2$.
Then, there are a finite subset $X'\subset\R^n$ and a
polynomial vector field $g:\R^n\rightarrow \R^n$ which has
asymptotically stable equilibria on the set $X$, 
saddles on the set $X'$, and no other equilibria, and moreover:
(1) every solution of $\dot x=g(x)$ converges to $X\bigcup X'$ and
(2) except for a measure-zero set of initial conditions, every solution
converges to an equilibrium in $X$.
\eps

\section{Proofs}

We prove Proposition~\ref{prop:main} by first treating the case of a
set $X\subset\R^2$ of the special form ${\cal X}\times \{0\}$, and then using a
coordinate change to reduce the general case to this one.  The proof of
Proposition~\ref{prop:weaker}, in contrast, only requires the coordinate
change, plus a trivial construction.

\subsection{A special case}

\bl{lemma:main1}
Let $\alpha :\R\rightarrow \R$ be a ${\cal C}^2$ function whose zeroes are all simple; that is, on
the set
\[
{\cal X} \, := \; \{x \st \alpha (x)=0\}
\]
it holds that $\alpha '(x)\not= 0$.
Introduce the following function $f:\R^2\rightarrow \R$:
\[
f(x,y) \, := \; 
 \left(\alpha (x) - \left[\alpha (x) - \alpha '(x)\right]^2y\right)^2 
   \,-\, \int \alpha (x)\left[\alpha (x) - \alpha '(x)\right]\,dx
\]
(where the last term denotes an arbitrary anti-derivative).
Consider the set of critical points of $f$,
\[
{\cal C}(f) \, := \; \left\{(x,y) \st f_x(x,y)=f_y(x,y) = 0 \right\} \,.
\]
Then:
\ben
\item
${\cal C}(f) = {\cal X} \times  \{0\}$.
\item
At each $(x,y)\in {\cal C}(f)$, the Hessian of $f$ is positive definite.
\een
As a consequence, $f$ has local minima at the points in ${\cal X} \times  \{0\}$, and no
other critical points.
\els

\bpr
For convenience, we introduce $\beta (x):=\alpha (x) - \alpha '(x)$, so that
\[
f(x,y) \; := \; 
 \left(\alpha (x) - \beta (x)^2y\right)^2 
   \,-\, \int \alpha (x)\beta (x)\,dx \,.
\]
We have:
\beqn
f_x(x,y) &=&
2 \left(\alpha (x) - \beta (x)^2y\right) \left(\alpha '(x) - 2\beta (x)\beta '(x)y\right)- \alpha (x)\beta (x)\\
f_y(x,y) &=&
- 2 \left(\alpha (x) - \beta (x)^2y\right) \beta (x)^2 \,.
\eeqn
Clearly, $\alpha (x)=0$ and $y=0$ imply $f_x(x,y)=f_y(x,y)=0$, so we must only
prove the converse.  Pick any $(x,y)\in {\cal C}(f)$.

From $f_y(x,y)=0$, we have that one of these must hold:
\be{eq:case1}
\beta (x)=0\,,
\ee
\be{eq:case2}
\beta (x)\not= 0 \;\mbox{and}\; \alpha (x) - \beta (x)^2y = 0\,.
\ee
If \eqref{eq:case1} holds, then $0=f_x(x,y)=2 \alpha (x)\alpha '(x)$ implies that
either $\alpha (x)=0$ or $\alpha '(x)=0$.
On the other hand,
the assumption of simple zeroes, ``$\alpha (x)=0 \,\Rightarrow \, \alpha '(x)\not= 0$''
can also be written in contrapositive form as
``$\alpha '(x)=0 \,\Rightarrow \,\alpha (x)\not= 0$,'' from which we have:
\beqn
\alpha (x)=0  & \Rightarrow  & \beta (x)=-\alpha '(x)\not= 0\\
\alpha '(x)=0 & \Rightarrow  & \beta (x)=\alpha (x)\not= 0 \,.
\eeqn
This rules out both $\alpha (x)=0$ and $\alpha '(x)=0$ when $\beta (x)=0$, and thus
case \eqref{eq:case1} cannot hold.

So case \eqref{eq:case2} holds, which means that $y=\alpha (x)/\beta (x)^2$.
On the other hand, when $\alpha (x) - \beta (x)^2y = 0$, we have that
$0 =f_x(x,y) = -\alpha (x)\beta (x)$, and since $\beta (x)\not= 0$, it follows that
$\alpha (x)=0$, from which it also follows that $y=0$.
We conclude that $(x,y)\in {\cal X}\times \{0\}$, as desired.

To prove positive definiteness of the Hessian on ${\cal C}(f)$, we must show that
on the set ${\cal C}(f)$, both $f_{yy}(x,y)>0$ and $\Delta (x,y)>0$,
where $\Delta  = f_{xx}f_{yy} - f_{xy}^2$.
In general:
\beqn
f_{xx} &=&  
2 \left(\alpha '(x) - 2\beta (x)\beta '(x)y\right)^2
+ 2 \left(\alpha (x) - \beta (x)^2y\right) \left(\alpha ''(x) - 2[\beta (x)\beta ''(x)+(\beta '(x))^2]y\right)\\
&& \quad\quad\quad - \alpha '(x)\beta (x)-\alpha (x)\beta '(x)\\
f_{yy} &=&
 2 \beta (x)^4\\
f_{xy} \;=\; f_{yx} &=&
- 2 \left(\alpha '(x) - 2\beta (x)\beta '(x)y\right)\beta (x)^2
- 4\left(\alpha (x) - \beta (x)^2y\right)\beta (x)\beta '(x)
\eeqn
so in particular, on the set ${\cal C}(f)$, since $\alpha (x)=0$ and $y=0$:
\beqn
f_{xx} &=&  
2 \left(\alpha '(x)\right)^2 - \alpha '(x)\beta (x)\\
f_{yy} &=&
 2 \beta (x)^4\\
f_{xy} = f_{yx} &=&
- 2 \left(\alpha '(x)\right)\beta (x)^2 \,.
\eeqn
Since on this set $\beta  = -\alpha '(x)$:
\beqn
f_{xx} &=& 3 \left[\alpha '(x)\right]^2\\
f_{yy} &=&  2 [\alpha '(x)]^4\\
f_{xy} = f_{yx} &=& - 2 \left[\alpha '(x)\right]^3
\eeqn
and thus
\[
\Delta (x,y) = 2 [\alpha '(x)]^6\,.
\]
As $\alpha '(x)\not= 0$ on the set ${\cal C}(f)$, it follows that both $f_{yy}>0$ and $\Delta >0$,
which completes the proof.
\epr

\subsection{A coordinate change}

\bl{lemma:main2}
Let $X$ be a finite subset of $\R^n$, $n\geq 2$.
Then, there exists a polynomial map $F:\R^n\rightarrow \R^n$ such that:
\be{eq:FmapsX}
F(X) \subseteq  \R\times \{0\}^{n-1} \,,
\ee
\be{eq:Finvertible}
F \;\; \mbox{has a polynomial inverse}.
\ee
\els

\bpr
For each pair of points $\xi \not= \eta $ in $X$, let 
$V_{\xi \eta }=\{p\in \R^n\st p^T(\xi -\eta )=0\}$.
Pick a point $p$ not belonging to the union of the finitely many hyperplanes
$V_{\xi \eta }$.
Choose any invertible mapping $T$ which has $p^T$ as its first row, and consider
the change of variables $z=Tx$.
Let $X= \{x^{(1)},\ldots ,x^{(k)}\}$ and $Z:=TX = \{z^{(1)},\ldots ,z^{(k)}\}$,
where $z^{(i)}=Tx^{(i)}$.
Since $x\mapsto p^Tx$ is one to one on the set $X$, the first coordinates of the
$z^{(i)}$'s are all distinct, that is
\[
z^{(i)} = (z^{(i)}_1,\ldots ,z^{(i)}_n)
\]
and $z^{(i)}_1\not= z^{(j)}_1$ for each $i\not= j$.
For each $j=2,\ldots ,n$, let $p_j$ be the Lagrange interpolation polynomial
that gives:
\[
p_j(z^{(i)}_1) = z^{(i)}_j\,,\quad\quad
i=1,\ldots ,k
\]
and define
\[
\Pi :\,\R^n\rightarrow \R^n:\,
\mypmatrix{z_1\cr z_2 \cr \vdots\cr z_n} \mapsto  
\mypmatrix{z_1\cr z_2 - p_2(z_1) \cr \vdots\cr z_n - p_n(z_1)} \,.
\]
This is invertible (with a polynomial inverse obtained by using
instead ``$z_j + p_j(z_1)$'' for each $j>1$).
Moreover, by definition, on the set $Z$ we have that the $j$th coordinate of
$P(z)$, $P(z^{(i)})_j = z^{(i)}_j - p_j(z^{(i)}_1) = 0$, $j>1$.
In other words, $P$ maps into $\R\times \{0\}^{n-1}$.
The proof is completed by picking the composition $F = \Pi \circ T$.
\epr

\subsection{Proof of Proposition~\protect{\ref{prop:main}}}

Let $X$ be a finite subset of $\R^n$, $n\geq 2$, pick $F$ as in
Lemma~\ref{lemma:main2}, and let ${\cal X}$ be such that $F(X)={\cal X}\times \{0\}^{n-1}$.
We will construct a polynomial $Q:\R^n\rightarrow \R$ whose only critical points are on
the set ${\cal X}\times \{0\}^{n-1}$, and the Hessian is positive definite there.
Then $P = Q\circ F$ will be as desired, because diffeomorphisms preserve critical
points and their signature.
(To be explicit: as $\nabla P(x) = \nabla Q(F(x))\cdot JF(x)$ and the Jacobian
$JF$ is everywhere nonsingular, critical points map to critical points.
Furthermore, at a critical point, the Hessian $H$ transforms as $J^T H J$,
so positive definiteness is preserved.)

Let $\alpha (x):=\prod_{i=1}^{k}(x-a_i)$, where ${\cal X}=\{a_1,\ldots ,a_k\}$.  This
$\alpha $ is as in Lemma~\ref{lemma:main1}; let $f$ be as there.
We define
\[
Q\left(x_1,\ldots ,x_n\right) := \;f(x_1,x_2) \,+\, \frac{1}{2}\sum_{i>2} x_i^2 \,.
\]
From $\partial Q/\partial x_i=0$ for $i>2$, the critical
points of $Q$ have $x_i=0$ for all $i>2$, so 
${\cal C}(Q)={\cal C}(f)\times \{0\}^{n-2} = {\cal X}\times \{0\}^{n-1}$.
Moreover, at these points, the Hessian of $Q$ is obtained by appending an
identity matrix to the Hessian of $f$, and thus it is also positive definite,
as required.
\qed

\subsection{Proof of Proposition~\protect{\ref{prop:weaker}}}

We use the same change of variables, so that, up to a diffeomorphism, we can
assume without loss of generality that
$X={\cal X}\times \{0\}^{n-1}$.
Let
\[
\gamma (x) :=\; - \prod_{i=1}^{k}(x-a_i) \prod_{i=1}^{k-1}(x-b_i) 
\]
where ${\cal X}=\{a_1 < \ldots  < a_k\}$ and the $b_i\in (a_i,a_{i+1})$ are arbitrary.
The scalar differential equation $\dot x=\gamma (x)$ has stable equilibria at
${\cal X}$, unstable at ${\cal X}'=\{b_1,\ldots ,b_{k-1}\}$, and no other equilibria.
We then define
\[
g(x) = (\gamma (x),-x_2,\ldots ,-x_n) \,.
\]
This satisfies the conclusions of the Lemma, with $X={\cal X}\times \{0\}^{n-1}$
and $X'={\cal X}'\times \{0\}^{n-1}$.
\qed

\edo
\begin{thebibliography}{AA}

\bibitem{robertson92}
Ian Robertson.
\newblock A polynomial with $n$ maxima and no other critical points.
\newblock Preprint; available at
http://www.mtholyoke.edu/\~{}adurfee/reu/92/reu92.htm, 1992.

\bibitem{durfee98}
Alan H. Durfee.
\newblock The index of grad $f(x, y)$.
\newblock \emph{Topology} {\bf 37}: 1339-1361, 1998

\end{thebibliography}
